# Nonparametric two-sample tests for increasing convex order

LUDWIG BARINGHAUS[*] and RUDOLF GRÜBEL[**]

*Institut für Mathematische Stochastik, Leibniz Universität Hannover, Postfach 60 09, D-30060 Hannover, Germany. E-mail:* [*]*lbaring@stochastik.uni-hannover.de;*
[**]*rgrubel@stochastik.uni-hannover.de*

Given two independent samples of non-negative random variables with unknown distribution functions $F$ and $G$, respectively, we introduce and discuss two tests for the hypothesis that $F$ is less than or equal to $G$ in increasing convex order. The test statistics are based on the empirical stop-loss transform, critical values are obtained by a bootstrap procedure. It turns out that for the resampling a size switching is necessary. We show that the resulting tests are consistent against all alternatives and that they are asymptotically of the given size $\alpha$. A specific feature of the problem is the behavior of the tests 'inside' the hypothesis, where $F \neq G$. We also investigate and compare this aspect for the two tests.

## 1. Introduction

One of the basic problems of actuarial mathematics and mathematical finance is the comparison of risks; see, for example, Kaas *et al.* (1994), Müller and Stoyan (2002) and Rolski *et al.* (1999). In order to introduce the aspect that we will deal with suppose that the random variables $X$ and $Y$ represent the loss associated with two insurance policies or other financial contracts; we assume throughout the paper that the random variables are non-negative and that they have finite mean. Clearly, if $X \leq Y$ then the comparison is a trivial task. However, the random quantities $X$ and $Y$ will in general not be directly comparable, and indeed, the comparison of risk is normally based on the respective distribution functions $F$ and $G$ of $X$ and $Y$.

A classical partial order for distributions is the *stochastic order*: We say that $F$ is less than or equal to $G$ in stochastic order and write $F \leq_{st} G$ (or, with some abuse of notation, $X \leq_{st} Y$) if

$$1 - F(x) \leq 1 - G(x) \qquad \text{for all } x \geq 0. \tag{1.1}$$

Many interesting and useful results are known for this concept, especially for its statistical aspects such as testing the hypothesis that $F \leq_{\text{st}} G$; see, for example, Chapter 6 in Conover (1971) or Kapitel 3,4 in Witting and Nölle (1970). However, in the context of risk comparison, stochastic order is too restrictive and it does not capture the important notion that risk should also depend on variability. For this and other reasons alternative notions of partial order for distributions have been investigated extensively, and the *increasing convex order* has turned out to play a major role in this application area. There are several equivalent definitions: We have $F \leq_{\text{icx}} G$ if

$$Ef(X) \leq Ef(Y) \qquad \text{for all } f \in \mathcal{F}_{\text{icx}}, \tag{1.2}$$

where $\mathcal{F}_{\text{icx}}$ denotes the set of increasing convex functions $f: \mathbb{R}_+ \to \mathbb{R}$, or if

$$E(X-t)^+ \leq E(Y-t)^+ \qquad \text{for all } t \geq 0, \tag{1.3}$$

see, for example, Theorem 1.5.7 in Müller and Stoyan (2002). Whereas (1.2) explains the terminology, (1.3) has the obvious interpretation in terms of reinsurance and stop-loss contracts. In fact, a straightforward application of Fubini's theorem leads to yet another condition. For this, we introduce the *stop-loss transform* $\text{SL}(F)$ or $F^{\text{SL}}$ associated with a distribution function $F$,

$$\text{SL}(F)(t) = F^{\text{SL}}(t) = \int_t^\infty (1 - F(x))\,\mathrm{d}x \qquad \text{for all } t \geq 0 \tag{1.4}$$

($F^{\text{SL}}$ is the integrated survival function; we prefer the notation $\text{SL}(F)$ whenever we want to emphasize the operator interpretation of the transform). Then the following is equivalent to (1.2) and (1.3):

$$F^{\text{SL}}(t) \leq G^{\text{SL}}(t) \qquad \text{for all } t \geq 0. \tag{1.5}$$

Indeed, the increasing convex order is often called *stop-loss order*, especially in an actuarial context. Obviously $X \leq_{\text{icx}} Y$ if $X \leq_{\text{st}} Y$ and $EY < \infty$; hence, for random variables with finite mean, increasing convex order is less restrictive than stochastic order. Further, by Jensen's inequality, we always have $X \leq_{\text{icx}} Y$ if $X \equiv EY$. This captures the fact that in an actuarial application, a fixed loss of magnitude $c$ would always be preferable to a random loss with mean $c$. In stochastic order distributions with the same mean are comparable only if they are identical.

When comparing actual risks an assertion such as $X \leq_{\text{icx}} Y$ may result from general considerations or may be the consequence of some model assumption. Here we take the view that data from previous contracts could serve as the basis for such a statement and we introduce two statistical tests that are applicable in this context. Such tests also have an obvious role in connection with model checking.

Formally, we assume that we have two independent samples $X_1, \ldots, X_m$ and $Y_1, \ldots, Y_n$, where the $X$-variables are independent and have distribution function $F$ and the $Y$-variables are independent and have distribution function $G$, and we consider the one-sided composite hypothesis $F \leq_{\text{icx}} G$ against the general alternative $F \not\leq_{\text{icx}} G$. In the special case where the distributions $F$ and $G$ are discrete and are concentrated on a fixed finite set the likelihood ratio test can be applied; this has been investigated in Liu and Wang (2003). Our approach in the general case is based on the plug-in principle. Let $1_A$ denote the indicator function associated with the set $A$. Regarding the empirical distribution functions

$$F_m = \frac{1}{m} \sum_{j=1}^{m} 1_{[X_j, \infty)} \quad \text{and} \quad G_n = \frac{1}{n} \sum_{k=1}^{n} 1_{[Y_k, \infty)},$$

of the two samples $X_1, \ldots, X_m$ and $Y_1, \ldots, Y_n$ as suitable estimators for $F$ and $G$, we are led to estimate $F^{\text{SL}}$ and $G^{\text{SL}}$ by the respective *empirical stop-loss transforms* $F_m^{\text{SL}}$ and $G_n^{\text{SL}}$. (We note in passing that, as SL is one-to-one, $F_m^{\text{SL}}$ is the nonparametric likelihood estimator for $F^{\text{SL}}$ in the same sense as $F_n$ has this property as an estimator for $F$; see, for example, Shorack and Wellner (1986), page 332f.) In view of the close relation of the present problem to the analogue for the classical stochastic order (see also Sections 3.1 and 3.5 below) natural candidates for test statistics are then the one-sided Kolmogorov–Smirnov test statistic for increasing convex order,

$$T_{m,n}^{\text{KS}} = \kappa_{m,n} \sup_{t \geq 0} \; (F_m^{\text{SL}}(t) - G_n^{\text{SL}}(t)), \tag{1.6}$$

and the one-sided Cramér–von Mises test statistic for increasing convex order,

$$T_{m,n}^{\text{CvM}} = \kappa_{m,n} \int_0^{\infty} (F_m^{\text{SL}}(t) - G_n^{\text{SL}}(t))^+ \, \mathrm{d}t, \tag{1.7}$$

where we have used the abbreviation

$$\kappa_{m,n} := \sqrt{\frac{mn}{m+n}} \tag{1.8}$$

for the multiplicative constants that are needed in order to obtain non-trivial limit distributions; see Section 2.2 below. In contrast to the classical stochastic order situation these statistics are not distribution-free, so a practical way to obtain critical values is needed. We propose to use the bootstrap; as we will see, this leads to a problem-specific variant that is a consequence of the nature of the hypothesis and that does not seem to arise in other situations where bootstrap tests have been suggested.

We have used the comparison of risks as a general motivation, but even within an actuarial context there are many different and specific applications of increasing convex order. An excellent source is the recent book by Denuit *et al.* (2005), where this order is used in connection with premium calculation, option pricing, the modeling of dependencies and the comparison of claim frequencies. Another important application area is renewal theory, where the stop loss transform of a lifetime distribution arises in connection with the stationary delay distribution and the limit distribution of the forward and backward recurrence times. Renewal theory in turn is a basic building block throughout the whole area of stochastic modeling.

In the next section we state our main results. Section 3 contains two examples and the results of a small simulation study; we also look at the practical calculation of the test statistics and we consider the relation to stochastic order in more detail. Proofs are collected in Section 4.

## 2. Main results

We first introduce some notation and list our general assumptions. Then we investigate the asymptotic behavior of the test statistics, both for situations where the hypothesis is true and situations where it is not. Our variant of obtaining bootstrap critical values is the topic of the next subsection. Finally, we combine these results into a description of the asymptotic behavior of the overall procedure, that is, where the test statistics are used together with the bootstrap critical values.

### 2.1. Preliminaries

In our main results a particular class of stochastic processes will be important: For a distribution function $H$ on $\mathbb{R}_+$, let $B_H = (B_H(t))_{t \geq 0}$ be a centered Gaussian process with covariance function $\rho_H$,

$$\rho_H(s,t) = \int_s^\infty \int_t^\infty (H(u \wedge v) - H(u)H(v))\, \mathrm{d}u\, \mathrm{d}v \qquad (2.1)$$

for all $s, t \geq 0$. Such a process can be obtained from a standard Brownian bridge $B = (B(t))_{0 \leq t \leq 1}$, which is a centered Gaussian process with covariance function

$$\operatorname{cov}(B(s), B(t)) = s \wedge t - s \cdot t,$$

via

$$B_H(t) = \int_t^\infty B(H(s))\, \mathrm{d}s \qquad \text{for all } t \geq 0. \qquad (2.2)$$

In particular, a process $B_H$ with the above properties exists if $\int x^2 H(\mathrm{d}x) < \infty$ (the role of moment conditions will be considered in more detail in the proofs section). The representation (2.2) could also informally be written as $B_H = \mathrm{SL}(B \circ H)$. In view of the facts that, first, $\sqrt{m}(F_m - F)$ converges in distribution to $B \circ F$ and, second, that SL is a linear operator it is not surprising that processes of this type appear in the context of distributional asymptotics for empirical stop-loss transforms.

The covariance structure of $B \circ H$ matches that of $(1_{[Z,\infty)}(t))_{t \geq 0}$, where $Z$ is a random variable with distribution function $H$. This (or a simple direct calculation) yields the alternative representation

$$\rho_H(s,t) = \operatorname{cov}((Z-s)^+, (Z-t)^+) \qquad \text{for all } s, t \geq 0 \qquad (2.3)$$

for the covariance function of $B_H$. In the proofs both (2.1) and (2.3) will be useful. We will also occasionally find it useful to extend these processes to the compactified half-line by $B_H(\infty) := 0$; a similar convention is used for the stop-loss transforms. Moment assumptions will imply that $B_H$ exists (as already mentioned above); that the paths are continuous (and hence bounded) functions on the compact set $K := [0, \infty]$, which is important in the Kolmogorov–Smirnov situation; or that these paths are integrable functions on $K$, which is necessary in the Cramér–von Mises case.

Let $F$ and $G$ be distribution functions on $\mathbb{R}_+$ with finite mean and let $(X_j)_{j \in \mathbb{N}}$ and $(Y_k)_{k \in \mathbb{N}}$ be two independent sequences of independent random variables, where the $X_j$'s have distribution function $F$ and the $Y_k$'s have distribution function $G$. Our convergence results refer to statistics $T_{m,n}$ that depend on the first $m$ of the $X$-variables and the first $n$ of the $Y$-values, where we generally assume that

$$\min\{m,n\} \to \infty, \qquad \frac{m}{m+n} \to \tau \in [0,1]. \tag{2.4}$$

We write $T_{m,n} \to_{\text{distr}} T$ if $T_{m,n}$ converges in distribution to $T$.

In order to describe the behavior of the tests inside the hypothesis we need two more definitions. Given $F$ and $G$,

$$A(F,G) := \{t \in [0,\infty] : F^{\text{SL}}(t) = G^{\text{SL}}(t)\} \tag{2.5}$$

denotes the set where the two stop-loss transforms are equal. Our general conventions imply that $A(F,G) \neq \varnothing$. Let $H := (1-\tau)F + \tau G$, with $\tau$ as in (2.4). Then

$$\gamma(H) := \sup\{t \in [0,\infty) : H(t) < 1\} \tag{2.6}$$

is the right end-point of the support of $H$. Finally, we generally assume that $F$ and $G$ are not degenerate and that the significance level $\alpha$ is less than $1/2$.

## 2.2. Asymptotic behavior of the test statistics

The first theorem in this subsection gives the limit distributions of $T_{m,n}^{\text{KS}}$ and $T_{m,n}^{\text{CvM}}$ in the case where the hypothesis is true.

**Theorem 1.** *Suppose that $F \leq_{\text{icx}} G$ and that $G$ has a finite second moment. Then, with $H = (1-\tau)F + \tau G$ and $B_H$ as in Section 2.1,*

$$T_{m,n}^{\text{KS}} \underset{\text{distr}}{\to} T^{\text{KS}} := \sup_{t \in A(F,G)} B_H(t).$$

*If $\int x^{4+\varepsilon} G(\mathrm{d}x) < \infty$ for some $\varepsilon > 0$, then*

$$T_{m,n}^{\text{CvM}} \underset{\text{distr}}{\to} T^{\text{CvM}} := \int_{A(F,G)} B_H^+(t)\,\mathrm{d}t.$$

On the alternative we have the following behavior.

**Theorem 2.** *Suppose that $F \not\leq_{\mathrm{icx}} G$. Then $T_{m,n}^{\mathrm{KS}}$ and $T_{m,n}^{\mathrm{CvM}}$ converge to $\infty$ with probability 1.*

### 2.3. Bootstrap critical values

The step from a test statistic $T$ to an actual statistical test $\phi = 1_{\{T>c\}}$ requires a critical value $c = c(\alpha)$ that depends on the chosen error bound $\alpha$ for the probability of wrongly rejecting the hypothesis. For stochastic order, due to the strong invariance properties and their consequences, the classical tests are level $\alpha$ tests even for finite sample sizes. In the present situation, however, we have to content ourselves with asymptotics, aiming for a procedure that satisfies the error bound asymptotically as the sample size(s) tend to infinity. At first, one might think of using Theorem 1 in the special case $F = G$ with the quantiles of the respective limit distribution as critical values. However, these limit distributions depend on the unknown $F$; also, it is not clear what the consequences are if we are inside the hypothesis in the sense that $F \leq_{\mathrm{icx}} G$, but $F \neq G$.

In this context resampling methodology provides a practicable solution. Recall that $F_m$ and $G_n$ are the empirical distribution functions associated with the samples $X_1, \ldots, X_m$ from $F$ and $Y_1, \ldots, Y_n$ from $G$. These random variables are the initial segments of two infinite sequences $(X_j)_{j \in \mathbb{N}}$ and $(Y_k)_{k \in \mathbb{N}}$ of random variables defined on some background probability space $(\Omega, \mathcal{A}, P)$; the almost sure statements below refer to $P$. Given the initial segments, let $\hat{Z}_{N,1}, \ldots, \hat{Z}_{N,N}$ be a sample of size $N := m + n$ from the (random) distribution function

$$H_{m,n} = \frac{n}{m+n} F_m + \frac{m}{m+n} G_n. \tag{2.7}$$

Note the 'size switching': In contrast to many two-sample bootstrap tests in the literature the basis for the resampling that we use here is a distribution that assigns higher probabilities to the values from the smaller sample (this aspect will be further discussed in Section 3.2). The following theorem shows that, with probability 1, the limit distribution for the respective test statistics is the same as in Theorem 1 in the special case $F = G$ (so that $A(F, G) = [0, \infty]$).

**Theorem 3.** *With $\hat{Z}_{N,1}, \ldots, \hat{Z}_{N,N}$ as above, let $\hat{F}_{N,m}$ and $\hat{G}_{N,n}$ be the empirical distribution functions associated with $\hat{Z}_{N,1}, \ldots, \hat{Z}_{N,m}$ and $\hat{Z}_{N,m+1}, \ldots, \hat{Z}_{N,N}$, respectively. Let*

$$\hat{T}_{m,n}^{\mathrm{KS}} := \kappa_{m,n} \sup_{t \geq 0} (\hat{F}_{N,m}^{\mathrm{SL}}(t) - \hat{G}_{N,n}^{\mathrm{SL}}(t)) \tag{2.8}$$

*and*

$$\hat{T}_{m,n}^{\mathrm{CvM}} := \kappa_{m,n} \int_0^\infty (\hat{F}_{N,m}^{\mathrm{SL}}(t) - \hat{G}_{N,n}^{\mathrm{SL}}(t))^+ \, \mathrm{d}t \tag{2.9}$$

be the bootstrap versions of the Kolmogorov–Smirnov and the Cramér–von Mises test statistics. Let $B_H$ be as in Section 2.1, with $H = (1-\tau)F + \tau G$. Then, with probability one,

$$\hat{T}_{m,n}^{\mathrm{KS}} \underset{\mathrm{distr}}{\to} \sup_{t\geq 0} B_H(t). \tag{2.10}$$

If $\int x^{4+\varepsilon} F(\mathrm{d}x) < \infty$ and $\int x^{4+\varepsilon} G(\mathrm{d}x) < \infty$ for some $\varepsilon > 0$, then, again with probability one,

$$\hat{T}_{m,n}^{\mathrm{CvM}} \underset{\mathrm{distr}}{\to} \int_0^\infty B_H^+(t)\,\mathrm{d}t. \tag{2.11}$$

The distributions of $\hat{T}_{m,n}^{\mathrm{KS}}$ and $\hat{T}_{m,n}^{\mathrm{CvM}}$ are functions of the data $X_1,\ldots,X_m, Y_1,\ldots,Y_n$. The familiar Monte Carlo procedure can be used to obtain approximations for the respective distribution functions and quantiles.

### 2.4. The bootstrap tests

Let $T^{\mathrm{KS}}$ and $T^{\mathrm{CvM}}$ be as in Theorem 1, and suppose now that $F = G$. In particular, the supremum and the integral are taken over the whole half-line. Let $c^{\mathrm{KS}}(F,\alpha)$ and $c^{\mathrm{CvM}}(F,\alpha)$ be the associated upper $\alpha$-quantiles, that is,

$$c^{\mathrm{KS}|CvM}(F,\alpha) := \inf\{q \geq 0 : P(T^{\mathrm{KS}|CvM} > q) \leq \alpha\}, \tag{2.12}$$

and let $\hat{c}_{m,n}^{\mathrm{KS}}(\alpha)$ and $\hat{c}_{m,n}^{\mathrm{CvM}}(\alpha)$ be the corresponding bootstrap estimates, that is, the upper $\alpha$-quantiles for $\hat{T}_{m,n}^{\mathrm{KS}}$ and $\hat{T}_{m,n}^{\mathrm{CvM}}$, respectively. We now consider the procedure that arises from using the test statistics defined in Section 1 together with these critical values. The following theorem is the main result of the paper. Our proof will give a slightly stronger result; we concentrate here on the statistically most relevant aspects. Let $\ell$ denote the Lebesgue measure.

**Theorem 4.** *Let $F$ and $G$ be distribution functions with finite second moment, put $S := A(F,G) \cap [0, \gamma(H))$. Then the limit*

$$\psi_{\mathrm{KS}}(\alpha; F, G) := \lim_{m,n\to\infty} P(T_{m,n}^{\mathrm{KS}} > \hat{c}_{m,n}^{\mathrm{KS}}(\alpha))$$

*exists, and*

$$\begin{aligned}
\psi_{\mathrm{KS}}(\alpha; F, G) &= \alpha && \text{if } F = G, \\
\psi_{\mathrm{KS}}(\alpha; F, G) &\leq \alpha && \text{if } F \leq_{\mathrm{icx}} G, \\
\psi_{\mathrm{KS}}(\alpha; F, G) &= 0 && \text{if and only if } F \leq_{\mathrm{icx}} G \text{ and } S = \varnothing, \\
\psi_{\mathrm{KS}}(\alpha; F, G) &= 1 && \text{if } F \not\leq_{\mathrm{icx}} G.
\end{aligned}$$

If $\int x^{4+\varepsilon} F(\mathrm{d}x) < \infty$ and $\int x^{4+\varepsilon} G(\mathrm{d}x) < \infty$ for some $\varepsilon > 0$, then

$$\psi_{\mathrm{CvM}}(\alpha; F, G) := \lim_{m,n \to \infty} P(T_{m,n}^{\mathrm{CvM}} > \hat{c}_{m,n}^{\mathrm{CvM}}(\alpha))$$

*exists, and*

$$\begin{aligned}
\psi_{\mathrm{CvM}}(\alpha; F, G) &= \alpha && \text{if } F = G, \\
\psi_{\mathrm{CvM}}(\alpha; F, G) &\leq \alpha && \text{if } F \leq_{\mathrm{icx}} G, \\
\psi_{\mathrm{CvM}}(\alpha; F, G) &= 0 && \text{if and only if } F \leq_{\mathrm{icx}} G \text{ and } \ell(S) = 0, \\
\psi_{\mathrm{CvM}}(\alpha; F, G) &= 1 && \text{if } F \not\leq_{\mathrm{icx}} G.
\end{aligned}$$

In words: The tests are asymptotically exact if the two distributions are the same; they are asymptotically of the preassigned level; they are consistent inside the hypothesis (where the meaning of 'inside' depends on the type of the test), and they are consistent against all alternatives. Whereas the Cramér–von Mises test requires a stronger moment assumption than the Kolmogorov–Smirnov test, it is asymptotically of level zero on a larger subset of the hypothesis. In particular, this 'consistency inside the hypothesis' can be used to distinguish between the two tests.

## 3. Examples, simulations and remarks

We first discuss the use of rank tests in Section 3.1. An example in Section 3.2 shows that the size switching in the resampling part mentioned after (2.7) is important; this example also illustrates the different behavior of the two tests inside the hypothesis. In Section 3.3 we provide two formulas that can be used to calculate the test statistics. Section 3.4 contains the result of a small simulation study. In the last subsection we show that, from an abstract statistical point of view, increasing convex order differs in an important way from the classical stochastic order.

### 3.1. Unsuitability of rank tests

What happens if we use the classical procedures designed for stochastic order, such as the Wilcoxon test and the Kolmogorov–Smirnov test, in connection with increasing convex order?

We write $\mathrm{Wei}(\beta)$ for the Weibull distribution with parameter $\beta > 0$ and $\mathrm{Exp}(\lambda)$ for the exponential distribution with parameter $\lambda > 0$; the respective distribution functions are $F(x) = 1 - \exp(-x^\beta)$ and $G(x) = 1 - \exp(-\lambda x)$, $x \geq 0$. The stop-loss transforms (integrated survival functions) of the Wei(2) and the Exp(1) distribution are

$$F^{\mathrm{SL}}(t) = \sqrt{\pi}(1 - \Phi(\sqrt{2}t)) \quad \text{and} \quad G^{\mathrm{SL}}(t) = \exp(-t), \qquad t \geq 0,$$

where $\Phi$ denotes the distribution function for the standard normal distribution. Figure 1 shows the two distribution functions on the left and the two stop-loss transforms on the right. Obviously, these distributions are not comparable in the usual stochastic order. However, we do have $F \leq_{\text{icx}} G$ in view of the Karlin–Novikov criterion (the mean of $F$ is less than or equal to the mean of $G$ and, for some $t_0$, $F(t) \leq G(t)$ for all $t \leq t_0$ and $F(t) \geq G(t)$ for all $t \geq t_0$; see, e.g., Theorem 3.2.4 in Rolski *et al.* (1999)).

Suppose we use the one-sided Wilcoxon test for testing $F \leq_{\text{icx}} G$ against the one-sided alternative $G \leq_{\text{icx}} F$, and the one-sided Kolmogorov–Smirnov test (in its version designed for stochastic order) for testing $F \leq_{\text{icx}} G$ against the general alternative $F \not\leq_{\text{icx}} G$. The first test rejects the hypothesis if $W_{m,n} > c_{m,n}$, where $W_{m,n} = \sum_{j=1}^{m} \sum_{k=1}^{n} 1_{\{X_j > Y_k\}}$ is the Wilcoxon–Mann–Whitney test statistic and $w_{m,n}$ is the $(1-\alpha)$-quantile of $W_{m,n}$ in the case where $F = G$. From $w_{m,n}/(mn) \to 1/2$ and

$$\frac{1}{mn}\sum_{j=1}^{m}\sum_{k=1}^{n} 1_{\{X_j > Y_k\}} \to \pi_{F,G}$$

in probability with $\pi_{F,G} := P(X_1 > Y_1) = \int (1 - F(x))G(\mathrm{d}x)$ it follows that $P(W_{m,n} > w_{m,n}) \to 1$ as $m, n \to \infty$ if $\pi_{F,G} > \frac{1}{2}$. Since

$$\pi_{\text{W}(2),\text{Exp}(1)} = \int_0^\infty \exp(-x - x^2)\,\mathrm{d}x > 0.54\ldots > \tfrac{1}{2},$$

we arrive at the conclusion that in the case $F = \text{Wei}(2)$ and $G = \text{Exp}(1)$ (where the hypothesis is true) the probability of an error of the first kind of the one-sided Wilcoxon test tends to 1 as $m, n \to \infty$.

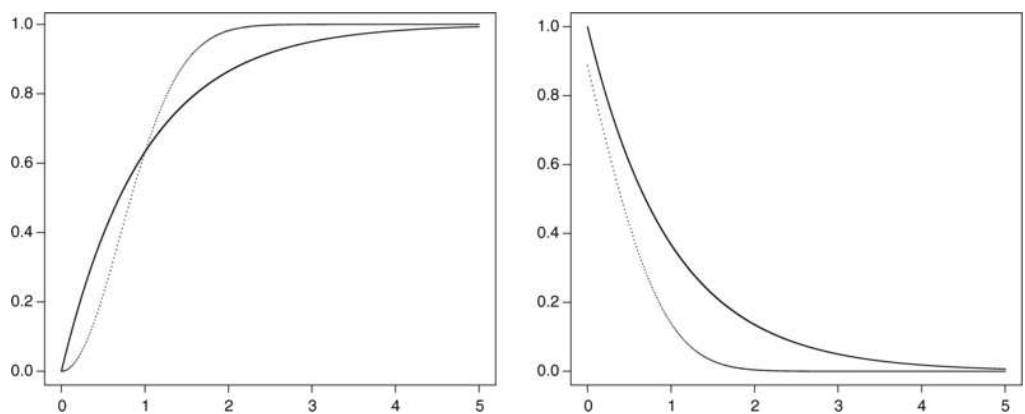

**Figure 1.** Distribution functions (left) and stop-loss transforms (right) for Exp(1) (solid line) and Wei(2) (dotted line).

The one-sided Kolmogorov–Smirnov test for stochastic order rejects the hypothesis if

$$\sup_{x\geq 0}(G_n(x) - F_n(x)) > k_{m,n},$$

where $k_{m,n}$ is the $(1-\alpha)$-quantile of $\sup_{x\geq 0}(G_n(x) - F_n(x))$ in the case where $F = G$. We always have $k_{m,n} \to 0$ and, with $F = \text{Wei}(2)$ and $G = \text{Exp}(1)$ again,

$$\sup_{x\geq 0}(G_n(x) - F_n(x)) \to \sup_{x\geq 0}(G(x) - F(x)) > 0$$

in probability as $m, n \to \infty$. Thus, for this test, too, the probability of an error of the first kind can become arbitrarily close to 1.

The two tests considered above are rank tests in the sense that the respective test statistics depend on the data only via the ranks of the sample variables in the pooled sample. We have seen that such tests may violate the bound $\alpha$ for the probability of wrongly rejecting the hypothesis $F \leq_{\text{icx}} G$ if used in conjunction with the critical values designed for stochastic order. Moreover, rank tests $\phi_{m,n}$ for increasing convex order that do respect the significance level $\alpha < 1$ (asymptotically) may fail to detect that the hypothesis is wrong, in the sense that they cannot be consistent against all alternatives. This follows from the fact that rank tests are invariant with respect to continuous, strictly increasing transformations $\Psi$ of the data, that is,

$$E_{F,G}\phi_{m,n} = E_{\tilde{F},\tilde{G}}\phi_{m,n}, \tag{3.1}$$

whenever $\tilde{F}$, $\tilde{G}$ are the distribution functions of $\Psi(X), \Psi(Y)$, respectively, and $X$ has distribution function $F$, $Y$ has distribution function $G$. Indeed, if $F$, $G$ and $\Psi$ are such that $F \leq_{\text{icx}} G$, $\tilde{F} \not\leq_{\text{icx}} \tilde{G}$, we have the (asymptotic) upper bound $\alpha$ on the left-hand side of (3.1) because of the bound on the probability of an error of the first kind, whereas consistency would require the limit 1 on the right-hand side as $m, n \to \infty$. An example for this situation can be obtained by choosing $F$ and $G$ to be the distribution functions for the uniform distribution on the intervals $(0, 1)$ and $(7/16, 9/16)$, respectively, together with $\Psi(x) = \sqrt{x}$. Then $F \leq_{\text{icx}} G$ follows with the Karlin–Novikov criterion, and a straightforward calculation shows that

$$\tilde{F}^{\text{SL}}(0) = \int x\tilde{F}(\mathrm{d}x) = \frac{9^{3/2} - 7^{3/2}}{12} > 2/3 = \int x\tilde{G}(\mathrm{d}x) = \tilde{G}^{\text{SL}}(0),$$

which implies that $\tilde{F} \not\leq_{\text{icx}} \tilde{G}$.

### 3.2. Size switching does matter

Suppose that the distributions of the $X$- and $Y$-variables are given by

$$P(X_j = 0) = P(X_j = 1) = \tfrac{1}{2} \quad \text{and} \quad P(Y_k = 0) = 1 - P(Y_k = 2) = \tfrac{3}{4},$$

respectively. The corresponding stop-loss transforms are

$$F^{\mathrm{SL}}(t) = \begin{cases} (1-t)/2, & 0 \leq t \leq 1, \\ 0, & t > 1, \end{cases} \quad \text{and} \quad G^{\mathrm{SL}}(t) = \begin{cases} (2-t)/4, & 0 \leq t \leq 2, \\ 0, & t > 2. \end{cases}$$

Obviously, $F \leq_{\mathrm{icx}} G$ and $A(F,G) = \{0\} \cup [2,\infty]$. In Section 2.3 we based the resampling on the distribution function

$$H_{m,n} := \frac{n}{m+n} F_m + \frac{m}{m+n} G_n.$$

We now investigate, for the above $F$ and $G$, what the consequences are if we resample from

$$H^0_{m,n} := \frac{m}{m+n} F_m + \frac{n}{m+m} G_n$$

instead. This corresponds to the standard resampling scheme where each of the items from the pooled sample is chosen with the same probability, so that the subsamples are represented proportional to their size when using $H^0_{m,n}$. This procedure has been suggested in Chapter 3.7.2 of van der Vaart and Wellner (1996) in connection with the general hypothesis $F = G$, for example, but of course, on the boundary $F = G$ of the hypothesis $F \leq_{\mathrm{icx}} G$ the size switching is asymptotically irrelevant.

Let $H = (1-\tau)F + \tau G$ and $H^0 := \tau F + (1-\tau)G$ be the respective limits of these distribution functions as $m, n \to \infty$, where we again assume that (2.4) holds. Let $B_H$ be as in Section 2.1 and let similarly $B_{H^0} = (B_{H^0}(t))_{t \geq 0}$ be a centered Gaussian process with covariance function

$$\rho_{H^0}(s,t) = \int_s^\infty \int_t^\infty H^0(u \wedge v) \, \mathrm{d}u \, \mathrm{d}v - \mathrm{SL}(H^0)(s)\mathrm{SL}(H^0)(t).$$

By Theorem 1 the test statistics $T^{\mathrm{KS}}_{m,n}$ converge in distribution to

$$T^{\mathrm{KS}} = \sup_{t \in A(F,G)} B_H(t) = B^+_H(0),$$

where we have used the fact that $B_H$ vanishes on $[2,\infty]$. By Theorem 3 the bootstrap estimators using size-switched resampling converge in distribution to

$$T^{\mathrm{switch}} = \sup_{t \geq 0} B_H(t)$$

with probability 1. The same arguments also work for standard resampling, leading to the distributional limit

$$T^{\mathrm{prop}} = \sup_{t \geq 0} B_{H^0}(t),$$

again with probability 1. Moreover, using the same arguments as in the proofs of Theorems 3 and 4 it can be shown that, with the $H^0$-based estimators $\hat{c}^{0,\mathrm{KS}}_{m,n}(\alpha)$ for the critical values,

$$\lim_{m,n\to\infty} P(T^{\text{KS}}_{m,n} > \hat{c}^{0,\text{KS}}_{m,n}(\alpha)) = P(T^{\text{KS}} > c^{\text{KS}}_0(\alpha)),$$

where $c^{\text{KS}}_0(\alpha)$ is the upper $\alpha$-quantile of the distribution of $T^{\text{prop}}$.

In the present simple case, where the distributions are concentrated on just two values, and using the fact that $B_H$ and $B_{H^0}$ can both be constructed from a standard Brownian bridge $B$, the limit processes essentially reduce to two-dimensional normal random variables. This enables us to work out some details. In fact, with the temporary notation $X := B(H(0)) = B((2+\tau)/4)$ and $Y := B(H(1)) = B((4-\tau)/4)$ we have that

$$B_H(t) = \int_t^\infty B(H(s))\,ds = \begin{cases} (1-t)X + Y, & 0 \le t < 1, \\ (2-t)Y, & 1 \le t < 2, \\ 0, & t \ge 2, \end{cases}$$

which gives $T^{\text{switch}} = X^+ + Y^+$. As $X$ and $Y$ are jointly normal with $\sigma_1^2 := \text{var}(X) = (4-\tau^2)/16$, $\sigma_2^2 := \text{var}(Y) = \tau(4-\tau)/16$ and $\text{cov}(X,Y) = \tau(2+\tau)/16$ we obtain, again with $\Phi$ the standard normal distribution function, and with $\rho := \text{cov}(X,Y)/(\sigma_1\sigma_2)$,

$$P(X^+ + Y^+ \le z) = \Phi\left(\frac{z}{\sigma_1}\right) - \int_{-\infty}^0 \left(1 - \Phi\left(\frac{z - \rho\sigma_2 t}{\sigma_2(1-\rho^2)^{1/2}}\right)\right)\Phi(dt)$$
$$- \int_0^{z/\sigma_1} \left(1 - \Phi\left(\frac{z - (\rho\sigma_2 + \sigma_1)t}{\sigma_2(1-\rho^2)^{1/2}}\right)\right)\Phi(dt)$$

for all $z \ge 0$. Similarly, replacing $\tau$ by $1 - \tau$, we obtain $T^{\text{prop}} = X_0^+ + Y_0^+$, now with $X_0 := B(H^0(0)) = B((3-\tau)/4)$ and $Y_0 := B(H^0(1)) = B((3+\tau)/4)$.

Table 1 gives the numerical results for some $\alpha$-values in the case $\tau = 3/4$ (the entries in columns 2–5 are rounded to four decimal points). The second and third columns contain the limiting critical values based on size-switched and proportional resampling, respectively, that is, the upper $\alpha$-quantiles $c^{\text{KS}}(\alpha)$ of the distribution of $T^{\text{switch}}$ and the upper $\alpha$-quantiles $c^{\text{KS}}_0(\alpha)$ of the distribution of $T^{\text{prop}}$. In columns 4 and 5 we have the corresponding asymptotic probabilities of rejection. The error bound is clearly violated for each of the $\alpha$-values in the proportional case. Due to the fact that $T^{\text{KS}} = B_H^+(0) = (X+Y)^+$ and $T^{\text{switch}} = X^+ + Y^+$ it is not surprising that in the case of size-switched resampling the limiting power is close to the given level for small $\alpha$: Conditionally on $X^+ + Y^+$ being large, we have $X, Y \ge 0$ and hence $(X+Y)^+ = X^+ + Y^+$ with high probability.

Turning to the behavior of the Cramér–von Mises test for increasing convex order we first note that $T^{\text{CvM}} = \int_{A(F,G)} B_H^+(t)\,dt = 0$. The limiting critical values are the upper $\alpha$-quantiles of $\int_0^\infty B_H^+(t)\,dt$ for size-switched resampling and $\int_0^\infty B_{H^0}^+(t)\,dt$ for proportional resampling. Since these quantiles are positive for each $\alpha \in (0, 1/2)$, the limiting power of the Cramér–von Mises test is 0 for both resampling schemes. With $F$ and $G$ as above we therefore also have an example for the different behavior of the two tests inside the hypothesis.

**Table 1.** Critical values and error probabilities

| $\alpha$ | $c^{\text{KS}}(\alpha)$ | $c_0^{\text{KS}}(\alpha)$ | $P(T^{\text{KS}} > c^{\text{KS}}(\alpha))$ | $P(T^{\text{KS}} > c_0^{\text{KS}}(\alpha))$ |
|---|---|---|---|---|
| 0.100 | 1.0134 | 0.8069 | 0.0999 | 0.1537 |
| 0.050 | 1.3004 | 1.0202 | 0.0500 | 0.0984 |
| 0.025 | 1.5495 | 1.2079 | 0.0250 | 0.0633 |

### 3.3. Alternative expressions for the test statistics

The definitions given in (1.6) and (1.7) result from applying the plug-in principle. From a practical point of view it is desirable to have alternative expressions that avoid the use of integrals. Suppose that the set

$$\{0, X_1, \ldots, X_m, Y_1, \ldots, Y_n\}$$

has $l$ different values and denote these, in increasing order, by $Z_1, \ldots, Z_l$.

**Proposition 1.** *Let $f_j := f(Z_j)$ for $j = 1, \ldots, l$, where*

$$f(t) := \frac{1}{m}\sum_{j=1}^{m}\max(X_j, t) - \frac{1}{n}\sum_{k=1}^{n}\max(Y_k, t),$$

*and let $J$ and $K$ be the complementary set of numbers $j = 1, \ldots, l-1$ such that $f_j \neq f_{j+1}$ and $f_j = f_{j+1}$, respectively. Then*

$$T_{m,n}^{\text{KS}} = \kappa_{m,n} \max_{1 \leq j \leq l-1} f_j^+$$

*and*

$$T_{m,n}^{\text{CvM}} = \kappa_{m,n}\left(\sum_{j \in J} \frac{1}{2}\frac{f_{j+1}^{+2} - f_j^{+2}}{f_{j+1} - f_j}(Z_{j+1} - Z_j) + \sum_{j \in K} f_j^+ (Z_{j+1} - Z_j)\right).$$

Clearly, these expressions can directly be used to implement the procedures on a computer.

### 3.4. Simulations

In order to obtain an impression of the performance of the tests we present the results of a simulation study. For each pair $F$ and $G$ of distributions chosen for the first and the second sample, respectively, an approximation of the power of the tests based on $T_{m,n}^{\text{KS}}$ and $T_{m,n}^{\text{CvM}}$ was obtained empirically by simulations with 1000 replications. For each of these samples we used 1000 resamples to get the corresponding approximate critical values. The test statistics were calculated with the formulas given in Section 3.3.

**Table 2.** Estimated powers of the tests; significance level $\alpha = 0.05$

| Distributions | | $m=50, n=30$ | | $m=30, n=50$ | | $m=50, n=50$ | |
|---|---|---|---|---|---|---|---|
| $F$ | $G$ | $T^{\mathrm{KS}}$ | $T^{\mathrm{CvM}}$ | $T^{\mathrm{KS}}$ | $T^{\mathrm{CvM}}$ | $T^{\mathrm{KS}}$ | $T^{\mathrm{CvM}}$ |
| unif$(0,1)$ | unif$(0,1)$ | 0.062 | 0.065 | 0.048 | 0.048 | 0.049 | 0.050 |
| unif$(\frac{1}{4},\frac{3}{4})$ | unif$(\frac{1}{4},\frac{3}{4})$ | 0.058 | 0.059 | 0.048 | 0.044 | 0.063 | 0.063 |
| Exp$(1)$ | Exp$(1)$ | 0.085 | 0.099 | 0.040 | 0.027 | 0.056 | 0.052 |
| Wei$(2)$ | Wei$(2)$ | 0.060 | 0.070 | 0.054 | 0.049 | 0.054 | 0.057 |
| $\Gamma(2)$ | $\Gamma(2)$ | 0.070 | 0.090 | 0.042 | 0.034 | 0.057 | 0.060 |
| $\Gamma(\frac{1}{2})$ | $\Gamma(\frac{1}{2})$ | 0.065 | 0.093 | 0.037 | 0.023 | 0.071 | 0.064 |
| unif$(\frac{1}{4},\frac{3}{4})$ | unif$(0,1)$ | 0.060 | 0.012 | 0.054 | 0.008 | 0.052 | 0.001 |
| Wei$(2)$ | Exp$(1)$ | 0.035 | 0.005 | 0.011 | 0.000 | 0.016 | 0.001 |
| unif$(0,2)$ | Exp$(1)$ | 0.090 | 0.044 | 0.068 | 0.011 | 0.078 | 0.013 |
| Exp$(4)$ | Par$(5)$ | 0.091 | 0.088 | 0.043 | 0.021 | 0.055 | 0.037 |
| unif$(0,1)$ | unif$(\frac{1}{4},\frac{3}{4})$ | 0.309 | 0.371 | 0.187 | 0.251 | 0.337 | 0.378 |
| Exp$(1)$ | Wei$(2)$ | 0.393 | 0.549 | 0.208 | 0.270 | 0.394 | 0.517 |
| Wei$(2)$ | Exp$(\frac{4}{3})$ | 0.284 | 0.117 | 0.272 | 0.061 | 0.346 | 0.087 |
| Exp$(\frac{4}{3})$ | Wei$(2)$ | 0.017 | 0.056 | 0.013 | 0.017 | 0.015 | 0.036 |
| $\Gamma(2)$ | Exp$(1)$ | 0.974 | 0.887 | 0.970 | 0.819 | 0.991 | 0.926 |
| Exp$(1)$ | $\Gamma(\frac{1}{2})$ | 0.812 | 0.642 | 0.806 | 0.493 | 0.892 | 0.634 |

Throughout, the sample sizes are $m, n \in \{30, 50\}$, and the significance level is $\alpha = 0.05$. The distributions $F$ and $G$ are chosen from the families of exponential distributions Exp$(\lambda)$, Weibull distributions Wei$(\beta)$, uniform distributions unif$(a, b)$ on finite intervals $(a, b)$, gamma distributions $\Gamma(\theta)$ with shape parameter $\theta > 0$ and scale parameter 1, and shifted Pareto distributions Par$(\eta)$ with parameter $\eta > 0$ and distribution function $1 - (1 + x)^{-\eta}$ for $x \geq 0$. Taking into account the relationships

$$\mathrm{unif}(1/4, 3/4) \leq_{\mathrm{icx}} \mathrm{unif}(0, 1) \leq_{\mathrm{st}} \mathrm{unif}(0, 2) \leq_{\mathrm{icx}} \mathrm{Exp}(1),$$

$$\Gamma(1/2) \leq_{\mathrm{st}} \mathrm{Exp}(1) \leq_{\mathrm{st}} \Gamma(2),$$

$$\mathrm{Wei}(2) \leq_{\mathrm{icx}} \mathrm{Exp}(1),$$

$$\mathrm{Exp}(4/3) \leq_{\mathrm{st}} \mathrm{Exp}(1),$$

$$\mathrm{Exp}(4) \leq_{\mathrm{icx}} \mathrm{Par}(5),$$

three groups of pairs of distributions are considered. Table 2 collects the results, where a numerical value such as 0.062 in the upper left corner means that 62 of the 1000 samples have led to a rejection; the entries therefore approximate the power of the respective test, here the Kolmogorov–Smirnov test of the hypothesis that $F \leq_{\mathrm{icx}} G$ in the case that the samples are both uniformly distributed on the unit interval.

The first group, lines 1–6 of the numerical values in Table 2, represents the boundary of the hypothesis, that is, cases where $F = G$. The two tests behave here essentially in the same way. The given level $\alpha = 0.05$ is exceeded slightly for sample sizes $m > n$ and the tests are nearly conservative for sample sizes $m < n$.

In the second group, lines 7–10, pairs from the interior of the hypothesis are considered, where $F \leq_{\text{icx}} G$ and $F \neq G$. For the first, third and fourth pairs of distributions it holds that $0 \in A(F, G)$ and $\int_{A(F,G)} B^+_{H_\tau}(t) \, dt = 0$. Thus, the limiting power of the Kolmogorov–Smirnov type test is positive and that of the Cramér–von Mises type test is 0. For the second pair of distributions $A(F, G)$ is empty and the limiting power of both tests is 0. Roughly, this behavior is reflected by the empirical power values. Results obtained for the special cases $F = \Gamma(1/2)$, $G = \text{Exp}(1)$ and $F = \text{Exp}(1)$, $G = \Gamma(2)$ are not shown in the table; here the empirical power value was found to be 0 for each combination of sample sizes.

The third group in lines 11–16 represents the alternative. None of the two tests turns out to be superior to the other one.

All in all, the simulations indicate that the asymptotic results in Section 2 lead to a feasible solution for test problems concerning increasing convex order even for relatively small sample sizes. From a practical point of view and especially with respect to implementation aspects, numerical stability, etc., our experience shows that the procedures are solid and uncomplicated. Of course, the resampling requires considerable computational effort, but with the current standard of hardware and software this is not really an issue. For example, we have done the entire simulation within the statistical environment R (see <http://cran.r-project.org>), which makes it possible to implement the procedures in less than 70 lines of code. For a specific data set with $n = m = 50$ and 1000 bootstrap resamples about 9 seconds were needed to calculate the test statistics and the critical value on a personal computer with an Intel Core 2 Duo processor running at 2.13 GHz. Of course, the time required could be shortened considerably by implementing the procedures directly in a programming language such as C.

## 3.5. Invariance considerations

One of the attractive properties of the classical tests for stochastic order, such as the Wilcoxon test and the Kolmogorov–Smirnov test mentioned in Section 3.1, is the fact that the distribution of the test statistics does not depend on $F$ if $F$ is continuous and $F = G$. As a result, exact critical values for finite sample sizes can be obtained. Also, stochastic order relates well to the ordering of real numbers: $X \leq_{\text{st}} Y$ implies the existence of a pair $X', Y'$ of random variables with $X' =_{\text{distr}} X$, $Y' =_{\text{distr}} Y$ and $X' \leq Y'$.

At the heart of these properties is the fact that the full group of all continuous and bijective (and hence strictly increasing) transformations $\Psi : \mathbb{R}_+ \to \mathbb{R}_+$ leaves the hypothesis $F \leq_{\text{st}} G$ invariant in the sense that

$$X \leq_{\text{st}} Y \iff \Psi(X) \leq_{\text{st}} \Psi(Y). \tag{3.2}$$

Such a strong property is not available in the present situation. Indeed, as the following theorem shows, only rescaling by positive constants leaves the hypothesis $F \leq_{\text{icx}} G$ invariant in the above sense.

**Theorem 5.** *Suppose that $\Psi:\mathbb{R}_+ \to \mathbb{R}_+$ is continuous and bijective, and such that*

$$X \leq_{\text{icx}} Y \quad \iff \quad \Psi(X) \leq_{\text{icx}} \Psi(Y) \tag{3.3}$$

*for all bounded random variables $X$ and $Y$. Then, for some $\alpha > 0$,*

$$\Psi(x) = \alpha x \qquad \text{for all } x \geq 0.$$

Hence the possibilities to reduce the test problem via invariance are limited. Of course, as we have already mentioned in Section 1, the increasing convex order has many properties that make it a suitable tool in connection with comparing risks; see also Chapter III in Kaas *et al.* (1994) or Section 1.5 in Müller and Stoyan (2002).

Whereas Theorem 5 may be regarded as a negative result as far as data reduction is concerned, it does say that the comparison of risks by increasing convex order will not be influenced by a change in monetary units, for example. It is straightforward to check that the tests in Section 2.4 are scale invariant in the sense that multiplying the data by a fixed positive constant does not affect the final decision, so our procedures respect the underlying symmetry of the problem.

### 3.6. Outlook

There are two aspects of our investigations that may be of interest in other situations, too. First, whereas consistency against alternatives is a property that is routinely investigated in connection with the asymptotic analysis of statistical tests, a similar property for the hypothesis is not usually considered. In the context of increasing convex order we have seen that this property can be used to distinguish between competing procedures, much in the same way as in a situation where there are different sets of alternatives on which two tests are consistent. Of course, efficiency considerations are then the logical next step; we plan to do this in a separate paper.

The other aspect, one that we have already pointed out repeatedly, is the fact that in order to obtain asymptotically correct critical values by resampling the resampling scheme has to be adapted to the hypothesis of interest – the standard procedure of generating artificial data by choosing values uniformly at random from the available observations may simply be wrong.

## 4. Proofs

Our main tool in connection with distributional asymptotics is the Hoffmann–Jørgensen theory of weak convergence, as described and developed in van der Vaart and Wellner (1996).

### 4.1. Proof of Theorem 1

For each $t \geq 0$, let $\phi_t : \mathbb{R}_+ \to \mathbb{R}$ be defined by $\phi_t(x) := (x - t)^+$, and let $\mathcal{F} := \{\phi_t : t \geq 0\}$. Then $\mathcal{F}$ is a Vapnik–Červonenkis class of index 2, see Lemma 2.6.16 in van der Vaart and Wellner (1996), and hence satisfies the uniform entropy condition. The function $\phi_0$ is an envelope for $\mathcal{F}$. Obviously, as in the transition from (1.3) to (1.5),

$$F_m^{\mathrm{SL}}(t) = \int \phi_t(x) F_m(\mathrm{d}x), \qquad G_n^{\mathrm{SL}}(t) = \int \phi_t(x) G_n(\mathrm{d}x) \qquad \text{for all } t \geq 0.$$

We now consider the empirical processes $U_m^F = (U_m(\phi))_{\phi \in \mathcal{F}}$ and $U_n^G = (U_n^G(\phi))_{\phi \in \mathcal{F}}$ associated with $X_1, \ldots, X_m$ and $Y_1, \ldots, Y_n$, respectively,

$$U_m^F(\phi) := \sqrt{m}\left( \int \phi \, \mathrm{d}F_m - \int \phi \, \mathrm{d}F \right), \qquad U_n^G(\phi) := \sqrt{n}\left( \int \phi \, \mathrm{d}G_n - \int \phi \, \mathrm{d}G \right)$$

for all $\phi \in \mathcal{F}$. The processes can be regarded as functions on the underlying probability space with values in $\ell_\infty(\mathcal{F})$, the Banach space of bounded real functions on $\mathcal{F}$ endowed with the supremum norm. We avoid the technical difficulties arising from the non-separability of $\ell_\infty(\mathcal{F})$ by noting that, for every distribution function $L$ with finite first moment, $t \mapsto \int \phi_t \, \mathrm{d}L$ is a bounded and continuous function on $[0, \infty)$, and that this function can be extended continuously to the compactified half-line $[0, \infty]$ by $\infty \mapsto 0$, again endowed with the supremum norm. Hence the empirical processes (and their distributional limits) all take values in a subspace of $\ell_\infty(\mathcal{F})$ that is a homeomorphic image of $C_0([0, \infty])$, the space of continuous functions $f : [0, \infty] \to \mathbb{R}$ with the property $f(\infty) = 0$. This subspace is closed and separable, and the empirical processes above are measurable with respect to the corresponding Borel $\sigma$-field. Indeed, we will freely switch between the representation of the 'time parameter' of the processes as an element $\phi$ of $\mathcal{F}$ or as an element $t$ of $K := [0, \infty]$.

By assumption, the envelope function is square integrable with respect to $F$ and $G$. Hence the function class $\mathcal{F}$ is $F$-Donsker and $G$-Donsker, so that

$$U_m^F \underset{\mathrm{distr}}{\to} B_F, \qquad U_n^G \underset{\mathrm{distr}}{\to} B_G \qquad \text{as } m, n \to \infty, \tag{4.1}$$

where $B_F = (B_F(\phi))_{\phi \in \mathcal{F}}$ and $B_G = (B_G(\phi))_{\phi \in \mathcal{F}}$ are centered Gaussian processes with covariance functions

$$\rho_F(\phi_s, \phi_t) = \int_s^\infty \int_t^\infty F(u \wedge v) \, \mathrm{d}u \, \mathrm{d}v - F^{\mathrm{SL}}(s) F^{\mathrm{SL}}(s),$$

$$\rho_G(\phi_s, \phi_t) = \int_s^\infty \int_t^\infty G(u \wedge v) \, \mathrm{d}u \, \mathrm{d}v - G^{\mathrm{SL}}(s) G^{\mathrm{SL}}(s),$$

respectively, for all $s, t \geq 0$; see Sections 2.1, 2.5 and 2.6 in van der Vaart and Wellner (1996). We now define the stochastic processes $V_{m,n}$, $m, n \in \mathbb{N}$, by

$$V_{m,n} := \kappa_{m,n}(F_m^{\mathrm{SL}} - G_n^{\mathrm{SL}})$$

$$= \sqrt{\frac{n}{m+n}} U_m^F - \sqrt{\frac{m}{m+n}} U_n^G + \kappa_{m,n}(F^{\mathrm{SL}} - G^{\mathrm{SL}}). \tag{4.2}$$

Note that, so far, the hypothesis $F \leq_{\mathrm{icx}} G$ has not been used. By the independence of the subsamples, and using (2.4), we obtain

$$\sqrt{\frac{n}{m+n}} U_m^F - \sqrt{\frac{m}{m+n}} U_n^G \underset{\mathrm{distr}}{\to} \sqrt{1-\tau} B_F - \sqrt{\tau} B_G =: \tilde{B},$$

where $\tilde{B}$ is again a centered Gaussian process, with covariance function

$$\tilde{\rho} = (1-\tau)\rho_F + \tau \rho_G.$$

Now let $B_H$ be as in Section 2.1. A straightforward calculation shows that

$$\rho_H(s,t) = \tilde{\rho}(s,t) + \tau(1-\tau)(F^{\mathrm{SL}}(s) - G^{\mathrm{SL}}(s))(F^{\mathrm{SL}}(t) - G^{\mathrm{SL}}(t)).$$

In particular, when restricted to the (compact, non-empty) set $A(F,G)$, the processes $B_H$ and $\tilde{B}$ have the same distribution.

We want to apply the continuous mapping theorem to obtain the limit distributions for

$$T_{m,n}^{\mathrm{KS}} = \sup_{t \geq 0} V_{m,n}(t), \qquad T_{m,n}^{\mathrm{CvM}} = \int_0^\infty V_{m,n}^+(t)\,\mathrm{d}t$$

in the case that $F \leq_{\mathrm{icx}} G$. Intuitively, it is clear that for the supremum or the positive part the values outside $A(F,G)$ are asymptotically irrelevant, as the last term in the basic equation (4.2) will tend to $-\infty$ for $t \notin A(F,G)$ in view of $\kappa_{m,n} \to \infty$. A formally correct argument for $T_{m,n}^{\mathrm{KS}}$ uses the Skorokhod almost sure representation together with the first part of the following auxiliary result.

**Lemma 1.** *Let $f_n$, $g_n$ ($n \in \mathbb{N}$), $g$, $h$ be continuous real functions on $K = [0,\infty]$ such that $f_n = g_n + c_n h$, where $(c_n)_{n \in \mathbb{N}}$ is a sequence of non-negative real numbers with $\lim_{n \to \infty} c_n = \infty$. Assume further that $h \leq 0$, that $A := \{h = 0\} \neq \varnothing$ and that $g_n$ converges uniformly to $g$. Then*

$$\lim_{n \to \infty} \sup_{t \in K} f_n(t) = \sup_{t \in A} g(t). \tag{4.3}$$

*and, for all $r < \infty$,*

$$\lim_{n \to \infty} \int_{[0,r]} f_n^+(t)\,\mathrm{d}t = \int_{[0,r] \cap A} g^+(t)\,\mathrm{d}t. \tag{4.4}$$

**Proof.** As the restriction of $g_n$ to the compact set $A$ converges uniformly to the restriction of $g$ to $A$ we have

$$\lim_{n \to \infty} \sup_{t \in A} g_n(t) = \sup_{t \in A} g(t).$$

Clearly,
$$\liminf_{n\to\infty} \sup_{t\in K} f_n(t) \geq \sup_{t\in A} g(t).$$

Suppose now that $(t_n)_{n\in\mathbb{N}} \subset K$ is such that $f_n(t_n) = \sup_{t\in K} f_n(t)$. We can find a subsequence $(n_k)_{k\in\mathbb{N}} \subset \mathbb{N}$ such that $t_{n_k} \to t_0 \in K$ and
$$\limsup_{k\to\infty} f_{n_k}(t_{n_k}) = \limsup_{n\to\infty} \sup_{t\in K} f_n(t).$$

A uniformly convergent sequence of continuous functions is equicontinuous, hence
$$\lim_{k\to\infty} g_{n_k}(t_{n_k}) = \lim_{k\to\infty} g_{n_k}(t_0) = g(t_0).$$

If $t_0 \notin A$ then $[t_0 - \varepsilon, t_0 + \varepsilon] \cap A = \varnothing$ for some $\varepsilon > 0$, and we would obtain $f_{n_k}(t_{n_k}) \to -\infty$ as $\sup_{t_0 - \varepsilon \leq t \leq t_0 + \varepsilon} h(t) < 0$ in view of the continuity of $h$. This shows that $t_0 \in A$, and hence
$$\sup_{t\in A} g(t) \geq g(t_0) = \limsup_{k\to\infty} g_{n_k}(t_{n_k})$$
$$\geq \limsup_{k\to\infty} f_{n_k}(t_{n_k}) = \limsup_{n\to\infty} \sup_{t\in K} f_n(t)$$
$$\geq \liminf_{n\to\infty} \sup_{t\in K} f_n(t) \geq \sup_{t\in A} g(t),$$

which proves (4.3).

For the proof of (4.4) we note that we have $f_n^+(t) = g_n^+(t) \to g^+(t)$ for $t \in A$ and that $f_n^+(t) \to 0$ for $t \notin A$. In view of the uniform convergence of $g_n$ to $g$ and $f_n^+ \leq g_n^+$ we have a constant upper bound for the sequence. As we integrate over a finite interval we can now use dominated convergence. □

For the proof of the Cramér–von Mises part of the theorem we temporarily abbreviate $A(F,G)$ to $A$. Using (4.4) together with an almost sure representation in the same way as we have used (4.3) for the proof of the Kolmogorov–Smirnov part we obtain
$$T^{\text{CvM}}_{m,n}(r) := \int_{[0,r]} V^+_{m,n}(t)\,dt \underset{\text{distr}}{\to} T^{\text{CvM}}(r) := \int_{A\cap[0,r]} B^+_H(t)\,dt$$

for all $r < \infty$. For the tails of the integrals we need another auxiliary result.

**Lemma 2.** *If $EX_1^{4+\varepsilon} < \infty$ for some $\varepsilon > 0$, then $\int_0^\infty \rho_F(t,t)^{1/2}\,dt < \infty$.*

**Proof.** Using (2.3) we obtain
$$\rho_F(t,t) = \text{var}((X-t)^+) \leq E((X-t)^+)^2$$
$$\leq \int_t^\infty x^2 F(dx) \leq \frac{1}{t^{2+\varepsilon}} EX^{4+\varepsilon}$$

for all $t > 0$. This shows that $\int_1^\infty \rho_F(t,t)^{1/2}\,dt < \infty$, and the assertion follows in view of the fact that $\rho_F$ is bounded (indeed, $\rho_F(t,t) \leq EX^2$). □

The moment assumption in the second part of the theorem and Lemma 2 together imply $\int_0^\infty \rho(t,t)^{1/2}\,dt < \infty$ for $\rho = \rho_F$, $\rho = \rho_G$ and $\rho = \rho_H$ (it is a straightforward consequence of (1.2) that finiteness of some moment of order $\gamma > 1$ for $G$ implies the same for $F$ if $F \leq_{\text{icx}} G$). It follows from

$$ET^{\text{CvM}} = \int_A EB_H^+(t)\,dt = \frac{1}{\sqrt{2\pi}}\int_A \rho_H(t,t)^{1/2}\,dt < \infty$$

that $T^{\text{CvM}}$ is finite with probability one. Clearly, $T^{\text{CvM}}(r) \to T^{\text{CvM}}$ almost surely as $r \to \infty$.

In order to obtain a uniform bound for the difference between $T_{m,n}^{\text{CvM}}(r)$ and $T_{m,n}^{\text{CvM}}$ we first note that $0 \leq V_{m,n}^+(t) \leq \tilde{V}_{m,n}^+(t)$, where

$$\tilde{V}_{m,n}^+(t) := \sqrt{\frac{n}{m+n}}U_m^F - \sqrt{\frac{m}{m+n}}U_n^G.$$

Clearly, $E\tilde{V}_{m,n}(t) = 0$ and $\text{var}(\tilde{V}_{m,n}(t)) \leq \rho_F(t,t) + \rho_G(t,t)$, and hence $EV_{m,n}^+(t) \leq (\rho_F(t,t) + \rho_G(t,t))^{1/2}$ for all $t \in A$. Suppose now that $\varepsilon > 0$ and $\delta > 0$ are given. We can then choose an $r > 0$ such that

$$\int_r^\infty (\rho_F(t,t) + \rho_G(t,t))^{1/2}\,dt < \varepsilon\delta,$$

and Markov's inequality gives

$$\limsup_{m,n\to\infty} P(|T_{m,n}^{\text{CvM}}(r) - T_{m,n}^{\text{CvM}}| > \varepsilon) \leq \frac{1}{\varepsilon}\int_r^\infty (\rho_F(t,t) + \rho_G(t,t))^{1/2}\,dt < \delta.$$

From this, the second assertion of the theorem now follows on using Theorem 4.2 in Billingsley (1968).

### 4.2. Proof of Theorem 2

The first two terms on the right-hand side of (4.2) are stochastically bounded. If $F \leq_{\text{icx}} G$ does not hold, then $F^{\text{SL}}(t) > G^{\text{SL}}(t)$ for all $t$ in some non-empty interval $(a,b) \subset \mathbb{R}$. As $\kappa_{m,n} \to \infty$, we therefore obtain

$$\sup_{t \geq 0} V_{m,n}(t) \to \infty, \qquad \int_0^\infty V_{m,n}^+(t)\,dt \to \infty$$

with probability one.

### 4.3. Proof of Theorem 3

We mimic the proof of Theorem 1 and give a somewhat condensed argument. With $\mathcal{F}$ as in Section 4.1, let $\hat{U}^F_{N,m} = (\hat{U}^F_{N,m}(\phi))_{\phi \in \mathcal{F}}$ and $\hat{U}^G_{N,n} = (\hat{U}^G_{N,n}(\phi))_{\phi \in \mathcal{F}}$,

$$\hat{U}^F_{N,m}(\phi) := \sqrt{m}\left(\int \phi \, d\hat{F}_{N,m} - \int \phi \, dF_m\right),$$

$$\hat{U}^G_{N,n}(\phi) := \sqrt{n}\left(\int \phi \, d\hat{G}_{N,n} - \int \phi \, dG_n\right),$$

be the empirical processes associated with the two parts of the resample. We need a statement analogous to (4.1), where now everything is conditional on the sequences $(X_j)_{j\in\mathbb{N}}$ and $(Y_k)_{k\in\mathbb{N}}$. Due to the size switching we cannot directly work with multiplier central limit theorems as, for example, in Chapter 3.6 of van der Vaart and Wellner (1996). Instead we base our proof on a central limit theorem for empirical processes where the base distribution varies with the sample size.

**Lemma 3.** *With probability one,*

$$\hat{U}^F_{N,m} \underset{\text{distr}}{\to} B_H, \qquad \hat{U}^G_{N,n} \underset{\text{distr}}{\to} B_H \qquad \text{as } m, n \to \infty.$$

**Proof.** It is obviously enough to prove the first part. Translated to the present situation the conditions (2.8.6) and (2.8.5) in Chapter 2.8.3 of van der Vaart and Wellner (1996) become

$$\limsup_{m,n\to\infty} E(\hat{Z}^2_{N,1} 1_{\{\hat{Z}_{N,1} > \varepsilon\sqrt{m}\}}) = 0 \qquad \text{for all } \varepsilon > 0$$

and

$$\lim_{m,n\to\infty} \sup_{s,t\geq 0} |\operatorname{var}(\phi_s(\hat{Z}_{N,1}) - \phi_t(\hat{Z}_{N,1}))^{1/2} - \operatorname{var}(\phi_s(Z_1) - \phi_t(Z_1))^{1/2}| = 0,$$

where $Z_1$ is a random variable with distribution $H$. Again this is to be interpreted in the sense that it should hold with probability one, conditionally on the sequences $(X_j)_{j\in\mathbb{N}}$ and $(Y_k)_{k\in\mathbb{N}}$.

For the first condition we can simply use the strong law of large numbers, together with the assumption that $F$ and $G$ have finite second moments.

The second condition will follow from

$$\lim_{m,n\to\infty} \sup_{t\geq 0} |E\phi_t(\hat{Z}_{N,1}) - E\phi_t(Z_1)| = 0 \tag{4.5}$$

and

$$\lim_{m,n\to\infty} \sup_{s,t\geq 0} |E\phi_t(\hat{Z}_{N,1})\phi_s(\hat{Z}_{N,1}) - E\phi_t(Z_1)\phi_s(Z_1)| = 0. \tag{4.6}$$

Obviously,

$$E\phi_t(\hat{Z}_{N,1}) - E\phi_t(Z_1) = \frac{n}{N}\frac{1}{m}\sum_{j=1}^{m}(X_j - t)^+ - (1-\tau)E(X_1 - t)^+$$

$$+ \frac{m}{N}\frac{1}{n}\sum_{k=1}^{n}(Y_k - t)^+ - \tau E(Y_1 - t)^+$$

and

$$E\phi_t(\hat{Z}_{N,1})\phi_s(\hat{Z}_{N,1}) - E\phi_t(Z_1)\phi_s(Z_1)$$

$$= \frac{n}{N}\frac{1}{m}\sum_{j=1}^{m}(X_j - t)^+(X_j - s)^+ - (1-\tau)E(X_1 - t)^+(X_1 - s)^+$$

$$+ \frac{m}{N}\frac{1}{n}\sum_{k=1}^{n}(Y_k - t)^+(Y_k - s)^+ - \tau E(Y_1 - t)^+(Y_1 - s)^+.$$

From this (4.5) and (4.6) follow as $\mathcal{F}$ and $\mathcal{F}_2 := \{\phi_s\phi_t : s, t \geq 0\}$ are easily seen to have finite bracketing numbers and thus are Glivenko–Cantelli classes; see Theorem 2.4.1 in van der Vaart and Wellner (1996). Now that we have checked the conditions the desired statement follows with Theorem 2.8.9 in van der Vaart and Wellner (1996). □

In analogy to (4.2) we now define the stochastic processes $\hat{V}_{N,m,n}$ by

$$\hat{V}_{N,m,n} := \sqrt{\frac{n}{m+n}}\hat{U}^F_{N,m} - \sqrt{\frac{m}{m+n}}\hat{U}^G_{N,n}.$$

By the conditional independence of the subsamples, and using (2.4) again, we obtain from Lemma 3 that $\hat{V}_{N,m,n} \to_{\text{distr}} B_H$ with probability one. Because of

$$\hat{T}^{\text{KS}}_{m,n} = \sup_{t\geq 0}\hat{V}_{N,m,n}(t)$$

the first part of the theorem now follows on using the continuous mapping theorem.

In the Cramér–von Mises case we start with a decomposition of the integral, as in the proof of Theorem 1. For the compact part we can again use the continuous mapping theorem. For the other part of the integral we need a bound for the conditional variances: We have

$$\text{var}(\hat{V}_{N,m,n}(t)) \leq \rho_{F_m}(t,t) + \rho_{G_n}(t,t),$$

and, as in the proof of Lemma 2,

$$\int_0^\infty \rho_{F_m}(t,t)^{1/2}\,dt \leq \frac{1}{m}\sum_{j=1}^{m}X_j^2 + \frac{1}{m}\sum_{j=1}^{m}X_j^{4+\varepsilon}\int_1^\infty t^{-1-(\varepsilon/2)}\,dt$$

and

$$\int_0^\infty \rho_{G_n}(t,t)^{1/2}\,\mathrm{d}t \leq \frac{1}{n}\sum_{k=1}^n Y_k^2 + \frac{1}{n}\sum_{k=1}^n Y_k^{4+\varepsilon}\int_1^\infty t^{-1-(\varepsilon/2)}\,\mathrm{d}t,$$

so we can use the strong law of large numbers to obtain an almost sure upper bound that does not depend on $m$ and $n$.

### 4.4. Proof of Theorem 4

We need the following properties of the limit distributions of the test statistics in the case that we are on the boundary of the hypothesis.

**Lemma 4.** *If $F = G$, with $F$ and $G$ as in Theorem 1, then the distribution functions of $T^{\mathrm{KS}}$ and $T^{\mathrm{CvM}}$ are strictly increasing on $[0,\infty)$ and continuous on $(0,\infty)$.*

**Proof.** Proposition 2 in Beran and Millar (1986) states that the distribution of the norm $\sup_{t\geq 0}|B_H(t)|$ has a density and a strictly increasing distribution function on $[0,\infty)$. Thus for each $a > 0$ we have

$$P\left(\sup_{t\geq 0} B_H(t) < a\right) \geq P\left(\sup_{t\geq 0}|B_H(t)| < a\right) > 0.$$

From Theorem 1 and the remark on page 854 in Tsirelson (1975) we deduce that the assertion of the proposition is true in the Kolmogorov–Smirnov case.

To obtain the corresponding assertion in the Cramér–von Mises case we argue as follows: Using Lemma 2 we obtain

$$E\int_0^\infty |B_H(t)|\,\mathrm{d}t = \int_0^\infty E|B_H(t)|\,\mathrm{d}t = \sqrt{\frac{2}{\pi}}\int_0^\infty \rho_H(t,t)^{1/2}\,\mathrm{d}t < \infty,$$

hence the process $B_H$ has integrable sample paths with probability one. Thus we can regard $B_H$ as a random element with values in a certain closed separable subspace of the Banach space $L^1([0,\infty),\mathrm{d}t)$. Again by Proposition 2 in Beran and Millar (1986) the distribution of the norm $\int_0^\infty |B_H(t)|\,\mathrm{d}t$ has a density and a strictly increasing distribution function on $[0,\infty)$. By the Hahn–Banach theorem there exists a sequence of functions $(e_n)_{n\in\mathbb{N}} \subset L^\infty([0,\infty),\mathrm{d}t)$ such that

$$\int_0^\infty |B_H(t)|\,\mathrm{d}t = \sup_{n\geq 1}\int_0^\infty B_H(t)e_n(t)\,\mathrm{d}t.$$

Noting that $B_H^+(t) = \frac{1}{2}(|B_H(t)| + B_H(t))$ we can then write

$$\int_0^\infty B_H^+(t)\,\mathrm{d}t = \tfrac{1}{2}\sup_{n\geq 1}\int_0^\infty B_H(t)(e_n(t)+1)\,\mathrm{d}t.$$

The random variables $\frac{1}{2}\int_0^\infty B_H(t)(e_n(t)+1)\,dt$, $n \in \mathbb{N}$, have centered normal distributions. Since

$$P\left(\sup_{n\geq 1} \frac{1}{2}\int_0^\infty B_H(t)(e_n(t)+1)\,dt < a\right) = P\left(\int_0^\infty B_H^+(t)\,dt < a\right)$$
$$\geq P\left(\int_0^\infty |B_H(t)|\,dt < a\right) > 0$$

for each $a > 0$ we can again apply the results in Tsirelson (1975) to see that the second assertion of the lemma is also true. □

Theorem 3 and the fact that the distribution function of the limit is strictly increasing together imply the convergence of the quantiles, that is,

$$\hat{c}_{m,n}^{\text{KS}}(\alpha) \to c^{\text{KS}}(H,\alpha), \qquad \hat{c}_{m,n}^{\text{CvM}}(\alpha) \to c^{\text{CvM}}(H,\alpha) \tag{4.7}$$

almost surely. From this, the last statement, the consistency on the full alternative, follows on using Theorem 2.

Suppose now that $F \leq_{\text{icx}} G$; let us temporarily write $F_{\text{KS}}$ and $F_{\text{CvM}}$ for the true limiting distribution functions of the test statistics and $\tilde{F}_{\text{KS}}$ and $\tilde{F}_{\text{CvM}}$ for the limiting distribution functions for the bootstrap estimators for the distributions of the test statistics. In both cases the limits can be obtained via the supremum or integral of the process $B_H$, over the whole of $\mathbb{R}_+$ for $\tilde{F}_{\text{KS}}$ and $\tilde{F}_{\text{CvM}}$ and over the subset $A(F,G)$ of $\mathbb{R}_+$ for $F_{\text{KS}}$ and $F_{\text{CvM}}$. Since $F$ and $G$ are assumed to be non-degenerate, we have that the supremum and the integral of the process $B_H$ over $\mathbb{R}_+$ are bounded from below by $B_H(0)$ and $\int_0^\infty B_H(t)\,dt$, which shows that $\tilde{F}_{\text{KS}}$ and $\tilde{F}_{\text{CvM}}$ are bounded from below in stochastic order by non-degenerate centered normal distributions. This implies that the quantiles $c^{\text{KS}}(H,\alpha)$ and $c^{\text{CvM}}(H,\alpha)$ are positive due to our general assumption that $\alpha < 1/2$. On the boundary of the hypothesis we have $H = F = G$, and $F_{\text{KS}} = \tilde{F}_{\text{CvM}}$ and $F_{\text{CvM}} = \tilde{F}_{\text{CvM}}$. Hence, using the continuity of $F_{\text{KS}}$ and $F_{\text{CvM}}$ on $(0,\infty)$

$$\lim_{m,n\to\infty} P(T_{m,n}^{\text{KS}} > \hat{c}_{m,n}^{\text{KS}}(\alpha)) = \alpha$$

and

$$\lim_{m,n\to\infty} P(T_{m,n}^{\text{CvM}} > \hat{c}_{m,n}^{\text{CvM}}(\alpha)) = \alpha,$$

which is the respective first statement in the theorem for the two tests.

For the proof of the remaining assertions we first note that $\text{var}(B_H(t)) = 0$ for $t \geq \gamma(H)$, so we may restrict the supremum or integral to $S = A(F,G) \cap [0,\gamma(H))$. The distribution functions $F_{\text{KS}}$ and $F_{\text{CvM}}$ have a jump of size 1 at 0 if $S = \varnothing$ or $\ell(S) = 0$ respectively, so that in this case

$$\lim_{m,n\to\infty} P(T_{m,n}^{\text{KS}} > \hat{c}_{m,n}^{\text{KS}}(\alpha)) = 0$$

and

$$\lim_{m,n\to\infty} P(T_{m,n}^{\mathrm{CvM}} > \hat{c}_{m,n}^{\mathrm{CvM}}(\alpha)) = 0.$$

On the other hand, if $S \neq \varnothing$ and $\ell(S) > 0$ we argue as in the proof of Lemma 4 to obtain that the distribution functions $F_{\mathrm{KS}}$ and $F_{\mathrm{CvM}}$ are strictly increasing on $[0,\infty)$ and continuous on the interval $(0,\infty)$. From this and $F_{\mathrm{KS}} \leq_{\mathrm{st}} \tilde{F}_{\mathrm{KS}}$ and $F_{\mathrm{CvM}} \leq_{\mathrm{st}} \tilde{F}_{\mathrm{CvM}}$ it follows with Theorem 1 that the limits

$$\psi_{\mathrm{KS}}(\alpha; F, G) = \lim_{m,n\to\infty} P(T_{m,n}^{\mathrm{KS}} > \hat{c}_{m,n}^{\mathrm{KS}}(\alpha)) = 1 - F_{\mathrm{KS}}(c^{\mathrm{KS}}(H,\alpha)) \leq \alpha$$

and

$$\psi_{\mathrm{CvM}}(\alpha; F, G) = \lim_{m,n\to\infty} P(T_{m,n}^{\mathrm{CvM}} > \hat{c}_{m,n}^{\mathrm{CvM}}(\alpha)) = 1 - F_{\mathrm{CvM}}(c^{\mathrm{CvM}}(H,\alpha)) \leq \alpha$$

exist and are positive.

### 4.5. Proof of Proposition 1

The function

$$f(t) := F_m^{\mathrm{SL}}(t) - G_n^{\mathrm{SL}}(t)$$
$$= \frac{1}{m}\sum_{j=1}^{m}(X_j - t)^+ - \frac{1}{n}\sum_{k=1}^{n}(Y_k - t)^+$$
$$= \frac{1}{m}\sum_{j=1}^{m}\max(X_j, t) - \frac{1}{n}\sum_{k=1}^{n}\max(Y_k, t)$$

is linear on the intervals $[Z_j, Z_{j+1})$, that is, of the form

$$f(t) = \alpha_j + \beta_j t, \qquad t \in [Z_j, Z_{j+1}),\ j = 1, \ldots, l-1,$$

with certain real $\alpha_j$ and $\beta_j$. Putting $f_j = f(Z_j)$, $j = 1, \ldots, l$, we have

$$T_{m,n}^{\mathrm{KS}} = \kappa_{m,n} \max_{1 \leq j \leq L-1} f_j^+.$$

Further, with $J$ and $K$ as in the statement of the proposition, we obtain the representation

$$T_{m,n}^{\mathrm{CvM}} = \kappa_{m,n}\left(\sum_{j \in J} \frac{1}{2} \frac{f_{j+1}^{+2} - f_j^{+2}}{f_{j+1} - f_j}(Z_{j+1} - Z_j) + \sum_{j \in K} f_j^+(Z_{j+1} - Z_j)\right).$$

### 4.6. Proof of Theorem 5

Let $Y$ be a bounded random variable and put $X \equiv EY$; then $X \leq_{\text{icx}} Y$ by Jensen's inequality. Hence if $\Psi$ satisfies (3.3) we must have $\Psi(EY) \leq E\Psi(Y)$ for all bounded random variables $Y$. From this property it easily follows, by taking $Y$ to be concentrated on two values, that $\Psi$ is convex. Now let $\Phi$ be the inverse of $\Psi$. Obviously, if (3.3) holds for $\Psi$, then it also holds for $\Phi$, so that $\Phi$ has to be convex, too. Now suppose that, for some $x_1, x_2 \in \mathbb{R}_+$, $\alpha \in (0,1)$,

$$\Psi(\alpha x_1 + (1-\alpha)x_2) < \alpha \Psi(x_1) + (1-\alpha)\Psi(x_2).$$

Then, with $y_i := \Psi(x_i)$, $i = 1, 2$, and using the fact that $\Phi$ is strictly increasing,

$$\alpha \Phi(y_1) + (1-\alpha)\Phi(y_2) < \Phi(\alpha y_1 + (1-\alpha)y_2),$$

which is a contradiction. Hence, $\Psi$ must satisfy

$$\Psi(\alpha x_1 + (1-\alpha)x_2) = \alpha \Psi(x_1) + (1-\alpha)\Psi(x_2)$$

for all $x_1, x_2 \in \mathbb{R}_+$ and all $\alpha \in (0,1)$. From the general assumptions on $\Psi$ it follows that $\Psi(0) = 0$. Taken together these properties imply that $\Psi(x) = x\Psi(1)$ for all $x \geq 0$.